\documentclass[11pt,twoside]{article}

\usepackage{amssymb}
\usepackage{latexsym}

\newtheorem{theorem}{Theorem}[section]
\newtheorem{lemma}[theorem]{Lemma}
\newtheorem{definition}[theorem]{Definition}
\newtheorem{proposition}[theorem]{Proposition}
\newtheorem{corollary}[theorem]{Corollary}
\newtheorem{notation}[theorem]{Notation}
\newtheorem{conjecture}{Conjecture}
\newtheorem{problem}{Problem}

\newenvironment{proof}{\paragraph{\it Proof.}}{$\square$\vskip0.4cm}
\newenvironment{remark}{\paragraph{\it Remark.}}{\vskip0.4cm}

\newcommand{\Qu}{{\mathbb Q}}
\newcommand{\M}{{\cal M}}
\newcommand{\N}{{\cal N}}
\newcommand{\Q}{{\cal Q}}
\newcommand{\J}{{\cal J}}
\newcommand{\C}{{\mathbb C}}
\newcommand{\R}{{\mathbb R}}

\newcommand{\Z}{{\mathbb Z}}
\newcommand{\E}{{\mathbb E}}
\newcommand{\K}{K_\Sigma}
\newcommand{\Ha}{{\cal H}}
\newcommand{\Hy}{{\mathbb H}}
\newcommand{\uPhi}{{\mathbf \Phi}}

\newcommand{\bM}{{\mathbf M}}
\newcommand{\calO}{{\cal O}}
\newcommand{\bL}{{\mathbb L}}
\newcommand{\Pic}{\mathop{\rm Pic}\nolimits}
\newcommand{\coker}{\mathop{\rm coker}\nolimits}
\newcommand{\im}{\mathop{\rm im}\nolimits}
\newcommand{\rank}{\mathop{\rm rank}\nolimits}
\newcommand{\ch}{\mathop{\rm ch}\nolimits}
\newcommand{\td}{\mathop{\rm td}\nolimits}
\newcommand{\Hom}{\mathop{\rm Hom}\nolimits}
\newcommand{\trace}{\mathop{\rm tr}\nolimits}
\newcommand{\End}{\mathop{\rm End}\nolimits}

\newcommand{\cE}{{\cal E}}
\newcommand{\compE}{{E\stackrel{\Phi}{\rightarrow}{E\otimes K}}}
\newcommand{\ucompE}{{\E_{\tM}\stackrel{\uPhi}{\rightarrow}
{\E_{\tM}\otimes\K}}}
\newcommand{\tM}{{\tilde{\M}}}
\newcommand{\tN}{{\tilde{\N}}}
\newcommand{\D}{{\mathbf D}}
\newcommand{\cL}{{\cal L}}

\begin{document}

%\title{Vanishing of intersection numbers on the moduli 
%space of Higgs bundles}
%\author{Tam\'as Hausel\\ School of Mathematics\\ Institute for
%Advanced Study\\ Olden Lane, Princeton, New Jersey 08540, USA
%\\ hausel@math.ias.edu}
%\maketitle

\thispagestyle{empty}

\leftline{\copyright~ 1998 International Press}
\leftline{Adv. Theor. Math. Phys. {\bf 2} (1998) 1011-1040}  

\vspace{0.4in}
\begin{center}
{\huge \bf Vanishing of intersection numbers }
\\[5mm]{\huge \bf on the moduli space of }
\\[5mm]{\huge \bf Higgs bundles}

\vspace{0.4in}

{\bf Tam\'as  Hausel}
\linebreak
School of Mathematics\\
Institute for Advanced Study\\ 
Olden Lane, Princeton\\
New Jersey 08540, USA\\
{\tt hausel@math.ias.edu} \\

\vspace{0.2in}
{\bf Abstract} \\
\vspace{0.1in}
\parbox[c]{4.5in}{\small \hspace{.2in} In this paper we consider the topological side of a
problem which is the analogue of Sen's S-duality testing
conjecture for Hitchin's moduli space
$\M$ of rank $2$ stable Higgs bundles of fixed determinant of odd
degree over a Riemann surface $\Sigma$. We prove that all 
intersection numbers in the compactly supported
cohomology of $\M$ vanish, i.e. ``there are no topological 
$L^2$ harmonic forms on $\M$''. This result
generalizes the well known vanishing of the Euler characteristic of the moduli
space of rank $2$ stable bundles $\N$ of fixed determinant of odd
degree over $\Sigma$. Our proof
shows that the vanishing of all intersection numbers of 
$H^*_{cpt}(\M)$ is given by
relations analogous to the Mumford relations in the cohomology ring of
$\N$.}
\end{center}
\renewcommand{\thefootnote}{}
\footnotetext{\small e-print archive: {\texttt http://xxx.lanl.gov/abs/math.AG/9805071}}
\renewcommand{\thefootnote}{\arabic{footnote}}

\section{Introduction}

Analyzing the conjectured S-duality in N=2 supersymmetric Yang-Mills theory, 
which is a proposed $SL(2,\Z)$ symmetry of the theory, Sen in \cite{sen} could
predict the dimension  of the space of $L^2$ harmonic forms $\Ha_k$ on 
the universal cover of the moduli space of magnetic monopoles of charge $k$,by speculating \newpage
 \pagenumbering{arabic}
\setcounter{page}{1012}

\pagestyle{myheadings}
\markboth{\it VANISHING OF INTERSECTION NUMBERS....}{\it T. HAUSEL }

\noindent that there must be an $SL(2,\Z)$ action on the space 
$\bigoplus\Ha_k$,
which represents bound electron states of the theory.

The moduli space of monopoles $M_k$ of charge $k$ is the space of
finite energy and charge $k$ solutions to the Bogomolny equations in 
$\R^3$, 
which can be interpreted as a reduction of the
self-dual $SU(2)$ Yang-Mills equations in $\R^4$. The space $M_k$ is a
non-compact manifold, with $\pi(M_k)=\Z_k$, and
has a natural hyperk\"ahler and complete metric on it, which comes from 
an abstract construction 
(the so-called hyperk\"ahler quotient construction,
cf. \cite{hitchin-et-al}) 
and known
explicitly only in the case $k=2$, when $M_2$ 
is called the Atiyah-Hitchin
manifold. (For further details see \cite{atiyah-hitchin}.) 

When $k=2$
Sen's conjecture says that $\dim(\Ha_2)=1$. By knowing the metric of
$M_2$ explicitly, Sen was able to
find a non-trivial $L^2$ harmonic form on the universal cover
$\tilde{M}_2$, giving some support
for his conjecture and in turn for S-duality.   

For higher $k$ 
Sen's conjecture says something about a metric which is not known 
explicitly. Nevertheless the statement is interesting from a 
mathematical point of view as the space of $L^2$ harmonic forms on a 
non-compact complete Riemannian manifold is not well understood. 

Hodge theory tells us that in the compact case the space of $L^2$ harmonic
forms is naturally isomorphic to the De-Rham cohomology of the manifold.
However in the non-compact case 
there is no such theory, and indeed 
the harmonic space depends crucially
on the metric.

Nevertheless some part
of Hodge theory survives for complete Riemannian
manifolds (cf. 
\cite{derham} Sect. 32 Theorem 24 and Sect. 35 Theorem 26),
such as the Hodge decomposition theorem which states that 
for a complete Riemannian manifold $M$ the space 
$\Omega^*_{L^2}$ of $L_2$ forms on $M$ has an orthogonal decomposition
$$\Omega^*_{L^2}=\overline{d(\Omega^*_{cpt})}\oplus \Ha^* \oplus 
\overline{\delta(\Omega_{cpt}^*)},$$ and also $\Ha^*=\ker(d)\cap
\ker(d^*)$. 
   
An easy corollary\footnote{Cf. \cite{segal-selby}} of these results 
says that the
composition  
$$H^*_{cpt}(M)\rightarrow \Ha^*\rightarrow
H^*(M)$$ is the forgetful map.

By calculating the image of $H^*_{cpt}(\tilde{M}_k)$ in $H^*(\tilde{M}_k)$ 
Segal and Selby could give a lower bound for the
harmonic forms on the moduli space of magnetic monopoles which coincides
with the dimension given by Sen's conjecture (see \cite{segal-selby}).
This purely mathematical result is thus a supporting evidence for 
the conjectured S-duality in $N=2$ SYM of theoretical Physics. 

In this paper we will investigate the analogue of Sen's conjecture for
Hitchin's moduli space $\M$ of rank $2$ Higgs bundles of fixed determinant of 
degree $1$ over a Riemann surface
$\Sigma$ of genus $g>1$. 
The space $\M$ is a simply connected non-compact manifold of
dimension
$12g-12$ with a
complete hyperk\"ahler metric on it, and was constructed by Hitchin
in \cite{hitchin1} by considering the solutions of the self-dual
Yang-Mills equations on $\R^4$ which are translation invariant in two
directions. 
Led by the similarities
between the spaces $M_k$ and $\M$ and their origin, we ask the
following question: 

\begin{problem} What are the $L^2$ harmonic forms on $\M$?
\label{l2}
\end{problem}
 
In this paper we prove the following:

\begin{theorem} The forgetful map\footnote{Unless otherwise stated
cohomology is meant with real coefficients.} 
$$j_\M:H_{cpt}^*(\M)\rightarrow H^*(\M)$$
is $0$. 
\label{main}
\end{theorem}

This says that unlike the case of 
$\tilde{M}_k$ the topology of $\M$ does not give
the existence of $L^2$ harmonic forms. We can state this fact
informally as: ``There are no {\em topological} $L^2$ harmonic forms on 
Hitchin's moduli space of Higgs bundles''. 

Segal and Selby's result together with 
Sen's conjecture suggest that for $\tilde{M}_k$ the topology gives all the
harmonic space. Led by this and supported by the discussion in
Subsection~\ref{toy}
we can formulate the following
conjecture: 

\begin{conjecture} There are no non-trivial $L^2$ harmonic forms on
Hitchin's moduli space of Higgs bundles. 
\label{conjecture}
\end{conjecture}

It would be interesting to see whether a physical argument could back
this conjecture. We know of one serious appearance 
of Hitchin's moduli space
of Higgs bundles in the Physics literature. In \cite{bershadsky-et-al}
a topological 
$\sigma$-model with target space $\M$ arises as certain limit of $N=4$
supersymmetric Yang-Mills theory. However it is not clear whether
$L^2$ harmonic forms on $\M$ have any physical interpretation in this theory.

Note that the conjecture does not hold for {\em parabolic} Higgs
bundles, as the toy example in Example 2 after Theorem 7.13 of 
\cite{hausel} shows. Note also that
Dodziuk's vanishing theorem \cite{dodziuk} shows that there are no
non-trivial 
$L^2$ {\em holomorphic} forms on $\M$, 
since the Ricci tensor of a hyperk\"ahler
metric is zero.  

From an algebraic geometrical point of view 
Theorem~\ref{main} can be interpreted as follows. First
of all it is really about middle dimensional cohomology,
because it is known that $\M$ does not have cohomology beyond the
middle dimension, and equivalently by Poincar\'e duality 
$\M$ does not have compactly supported cohomology below the middle dimension. 
Thus the main content of Theorem~\ref{main} is the vanishing of
the canonical map $j_\M:H_{cpt}^{6g-6}(\M)\rightarrow H^{6g-6}(\M)$
between $g$ dimensional spaces (cf. Corollary 5.4 in \cite{hausel}). 
This in turn is
equivalent to the vanishing of the intersection form on
$H_{cpt}^{6g-6}(\M)$, i.e. to the vanishing of $g^2$ 
intersection numbers. 

There are $g+1$ intersection numbers whose 
vanishing follows easily.
One vanishing is obtained by recalling that the moduli space $\N$ of stable
bundles of real dimension $6g-6$ 
sits inside $\M$ with normal bundle $T^*_\N$, thus its
self-intersection number is its Euler characteristic up to sign, which is
known to vanish. 
 
The other $g$ vanishings follow from the fact that the ordinary
cohomology class of the Prym variety, the generic fibre of the Hitchin
map, is $0$, i.e. it is in the kernel of $j_\M$. This can be seen
by thinking of the Hitchin map as a section of the trivial
rank $3g-3$ vector bundle on $\M$ and considering the ordinary cohomology class
of the Prym variety as the Euler class of this trivial vector bundle, and as
such, the ordinary cohomology class of the Prym variety is trivial
indeed.   
Note that for the case $g=2$, the above vanishings are already enough
to have $j_\M=0$. (Cf. Example 2 after Theorem 7.13 in \cite{hausel}.)
 
The vanishing of the rest of the $g^2$ intersection numbers
on $\M$, proved in the present paper for any genus, 
can be considered as a generalization of these facts.  

 The structure of this paper is as follows: In the next section we
describe
the cohomology of certain moduli spaces. In Section~\ref{hyper} we
define hypercohomology groups and related notions. 
Then in Section~\ref{vanishin} we develop the theory of stable Higgs
bundles analogously to the stable vector bundle
case, and prove an important vanishing theorem. In
Section~\ref{univbundle}
we prove that $\tM$ is a fine moduli space, and define certain
universal bundles. In Section~\ref{virdirac} we construct the
virtual Dirac bundle as the analogue of the virtual
Mumford bundle,  and show that it can be considered as the degeneracy
sheaf of a homomorphism of vector bundles. In
Section~\ref{degloc} we determine the degeneracy locus of the
above homomorphism in terms of the components of the nilpotent
cone. Finally in Section~\ref{proof} we prove our main
Theorem~\ref{main} using Porteous' formula for the degeneracy locus of
the virtual Dirac bundle.

\paragraph{Acknowledgements.} First of all I would like to thank my
supervisor Nigel Hitchin for suggesting Problem~\ref{l2},
and for his help and encouragement. I am grateful to Michael
Thaddeus for his inspiring paper \cite{thaddeus1}, 
enlightening communications and
his constant interest in my work. I am also indebted to Manfred Lehn for
the idea of the proof of Theorem~\ref{deg}. I have found conversations
with Michael Atiyah, Frances Kirwan and Graeme Segal very stimulating. 
I thank the Mathematical Institute
and St. Catherine's College, Oxford for their hospitality during the
preparation of this work. Finally I thank Trinity College, Cambridge
for financial support.

\section{Moduli spaces and their cohomology}
\label{cohmod}

The central object of this paper is a  
fixed, smooth and complex projective curve $\Sigma$ of genus $g\geq 2$. We
also fix a point $p\in \Sigma$.

An additive basis of $H^*(\Sigma)$: $1\in H^0(\Sigma)$, 
$e_i\in H^1(\Sigma),\  i=1,..,2g$ and the fundamental cohomology class 
$\sigma\in H^2(\Sigma)$ with the properties that 
$e_i\wedge e_{i+g} = -e_{i+g}\wedge e_i =\sigma$ for $i=1,..,g$ and otherwise
$e_i\wedge e_j=0$, will be fixed throughout this paper.

\subsection{The Jacobian $\J$}

The moduli space of 
line bundles of degree $k$ over $\Sigma$ is the Jacobian $\J_k$. This is an 
Abelian variety of dimension $g$. Tensoring by a fixed line bundle of degree
$k-l$ gives an isomorphism between $\J_l$ and $\J_k$.
We will write $\J$ for $\J_1$. 

Being a torus 
$H^*(\J_k)$ is a free exterior algebra on $2g$ classes $\tau_i\in H^1(\J_k)$
defined by the formula $$c_1(\bL_k)=k\otimes\sigma + 
\sum_{i=1}^{2g}\tau_i\otimes e_i \in H^2(\J_k\times \Sigma)\cong
\sum_{r=0}^{2} H^r(\J_k)\otimes H^{2-r}(\Sigma).$$ Here $\bL_k$ is
the normalized Poincar\'e bundle, or universal line bundle over 
$\J_k\times \Sigma$. Universal means that for any $L\in \J_k$: 
$$\bL_k\mid_{\{L\}\times \Sigma}\cong L$$
and normalized means that $\bL_k\mid_{\J_k\times \{p\}}$ is trivial 
(cf. \cite{arbarello-et-al}).

\subsection{Moduli space of Abelian Higgs bundles $T^*_\J$}
\label{toy}

As a toy example for the discussions in the Introduction, we
consider here the moduli space of Abelian Higgs bundles. 

The tangent bundle of $\J$ is canonically isomorphic to $\J\times 
H^1(\Sigma,\calO_\Sigma)$. Thus by Serre duality 
$T^*_{\J}\cong \J\times
H^0(\Sigma,K)$ canonically. An element $\Phi\in (T^*_\J)_L\cong
H^0(\Sigma,K)$, can be thought of as a rank $1$ Higgs bundle:
$\cL=L\stackrel{\Phi}{\rightarrow} L\otimes K$
(cf. Definition~\ref{higgs}). 
Thus
we can think of $T^*_\J$ as the moduli space of rank $1$ Higgs
bundles. 

The cohomology of $T^*_\J$ is isomorphic to that of $\J$. However there
is an extra piece of cohomological information namely the intersection
numbers in the compactly supported cohomology or in other words the
map:
$$j_\J:H^*_{cpt}(T^*_\J)\rightarrow H^*(T^*_\J).$$ 
Clearly this map is interesting only in the middle dimension, where 
both $H^{2g}_{cpt}(T^*_\J)$ and $H^{2g}(T^*_\J)$ are
one-dimensional. However the Euler characteristic of $\J$ is clearly
$0$, thus the self-intersection number of the zero section of $T^*_\J$ is
$0$, which shows that $j_\J$ vanishes. 

Consider the Riemann 
metric on $T^*_\J\cong  \J\times H^0(\Sigma,K)$ which is the product 
of the flat metrics on the two terms (this is the metric which we get
if we perform Hitchin's work in \cite{hitchin1} for the Abelian case).
From the $L^2$-vanishing theorem of Dodziuk \cite{dodziuk}, since the
metric is flat 
there are no non-trivial $L^2$ harmonic forms on $T^*_\J$, thus
in the Abelian Higgs case the topology gives the harmonic space, as 
conjectured for the rank $2$ Higgs moduli space in 
Conjecture~\ref{conjecture}
and for the universal cover of the 
moduli space of magnetic monopoles in \cite{sen}.

\subsection{The moduli space of rank $2$ stable bundles $\N$}

We denote by $\tilde{\N}$ the fine moduli space of rank $2$ 
stable bundles with degree $1$ over $\Sigma$. 
It is a smooth projective variety of dimension
$4g-3$. The determinant gives a map 
$det_\N:\tilde{\N}\rightarrow \J$. For any 
$\Lambda\in \J$ the fibre $det_\N^{-1}(\Lambda)$ will be
denoted by $\N_{\Lambda}$, which is a smooth projective variety
of dimension $3g-3$. The map $f:\N_{\Lambda_1}\rightarrow
\N_{\Lambda_2}$ given by $f(E)= E\otimes (\Lambda_2\otimes
\Lambda_1^*)^{1/2}$, where $(\Lambda_2\otimes \Lambda_1^*)^{1/2}$ is a
fixed square root of $\Lambda_2\otimes \Lambda_1^*$, is an isomorphism
between $\N_{\Lambda_1}$ and $\N_{\Lambda_2}$. Hence we will write
$\N$ for $\N_\Lambda$, when we do not want to emphasize the fixed line bundle
$\Lambda$. 

Let $G:=H^1(\Sigma,\Z_2)\cong \Z^{2g}\cong \ker(\sigma_2)$, 
where $\sigma_2:\J_0\rightarrow \J_0$ is given by $\sigma_2(L)=L^2$. Now
$G$ acts on $\N$ and $\J$ by tensoring with the corresponding line bundle in 
$\ker(\sigma_2)$ and also on $\N\times \J$ by the diagonal action. 
Then (cf. (9.5) of \cite{atiyah-bott}) we have 
\begin{eqnarray}
\tN=(\N\times \J)/G.\label{split}
\end{eqnarray} 
Because 
$G$ acts trivially on $H^*(\J)$ and on $H^*(\N)$ (the latter was
first proved in \cite{harder-narasimhan}) we see that as rings 
\begin{eqnarray}
H^*(\tN)\cong \left( H^*(\N) \otimes H^*(\J) \right) ^G\cong H^*(\N)\otimes
H^*(\J).
\label{splitN}\end{eqnarray}
Thus for understanding the cohomology ring $H^*(\tN)$ it is enough
to know the cohomology ring $H^*(\N)$. The latter is 
multiplicatively generated by classes
$\alpha_\N\in H^2(\N)$, $\psi^i_\N\in H^3(\N)$ and $\beta_\N\in
H^4(\N)$, which appear in the K\"unneth decomposition of $c_2(\End(\E_\N))$:
\begin{eqnarray} 
c_2(\End(\E_\N))=2\alpha_\N\otimes\sigma+
\sum_{i=1}^{2g}4\psi_\N^i\otimes e_i -
\beta_\N\otimes 1 \label{kunneth}
\end{eqnarray}
in $H^4(\N\times \Sigma)\cong
\sum_{r=0}^{4} H^r(\N)\otimes H^{4-r}(\Sigma)$.

Here $\E_\N$ is the normalized 
rank $2$ universal bundle over $\N\times \Sigma$, i.e.
$c_1(\E_\N)=\alpha_\N$ and $\E_\N\mid_{\{E\}\times \Sigma}\cong E$ for every $E \in \N$.

The ring $H^*(\N)$ is described in terms of the so called 
{\em Mumford relations}. To explain this
consider the {\em virtual Mumford bundle} 
\begin{eqnarray}
\bM=-{\pi_\tN}_!(\E_\tN\otimes \pi_\Sigma^*(L_p^{-1}))&=&
-R^0{\pi_\tN}_*(\E_\tN\otimes \pi_\Sigma^*(L_p^{-1}))\nonumber\\
&& +R^1{\pi_\tN}_*
(\E_\tN\otimes 
\pi_\Sigma^*(L_p^{-1}))\in K(\tN).\nonumber
\end{eqnarray}
 Using standard properties of stable bundles it can be shown that
$R^0$ vanishes. Thus $\bM$ is a vector bundle of rank $2g-1$. Its total Chern
class is a complicated\footnote{It was calculated by Zagier in
\cite{zagier}.} 
polynomial of the universal classes.
Since $\rank(\bM)=2g-1$, the Chern class
$c_{2g+r}(\bM)\in H^{4g+2r}(\tN)$ vanishes for $r\geq 0$. According to
(\ref{splitN}), the cohomology of $\tN$ is the tensor product of
$H^*(\J)$ and $H^*(\N)$. Thus if we write $\tau_S=\prod_{i\in
S}\tau_i\in H^{|S|}(\J)$
for $S\subset \{1\dots 2g\}$ and 
$$c_{2g+r}(\bM)=\sum_{S\subset \{1\dots 2g\}} \zeta^r_S\otimes \tau_S$$ in the
K\"unneth decomposition of (\ref{splitN}) then we get the vanishing of
each $\zeta^r_S$. Thus for every $r\geq 0$ and $S\subset \{1\dots 2g\}$ 
we get a relation
\begin{eqnarray}\zeta^r_S\in \Qu[\alpha,\beta,\psi_i]\label{mumfordrelation}\end{eqnarray} 
of degree $4g+2r-|S|$. The polynomials 
$\zeta^r_S$ are called the
{\em Mumford relations}.
 
Mumford conjectured, and it was first
proved by Kirwan in \cite{kirwan}, that the Mumford relations constitute
a complete set of relations of the cohomology ring of $\N$.

By now a complete
description of the Mumford relations and the ring structure of $H^*(\N)$ 
is available (see
\cite{baranovsky},
\cite{king-newstead}, \cite{siebert-tian} and \cite{zagier} 
and also \cite{thaddeus2} for an introduction to the topology of $\N$).

\subsection{The moduli space of rank $2$ stable Higgs bundles $\M$}

We denote by $\tM_{2k-1}$ 
the coarse moduli space of rank $2$ stable Higgs bundles\footnote{For 
definitions see Section~\ref{vanishin}.} 
with degree $2k-1$ over $\Sigma$, 
which was constructed as $\M(2,2k-1,K)$ in
\cite{nitsure}. For a fixed $\Sigma$ they are all isomorphic to each
other. We write $\tM$ for $\tM_1$.
It is a smooth, non-projective,  
quasi-projective variety of dimension
$8g-6$. 

The determinant gives a map 
$det_\M:\tM\rightarrow T^*_\J$, defined by 
$det_\M(E,\Phi)=(\Lambda^2E,\trace(\Phi))$. For any 
$\cL\in T^*_\J$ the fibre $det_\M^{-1}(\cL)$ will be
denoted by $\M_{\cL}$. Just as in the stable vector bundle 
case any two fibres of
$det_\M$
are 
isomorphic. Usually we will write $\M$ for
$\M_\cL$, when the Abelian Higgs bundle $\cL$ has zero Higgs field. 

Our main concern in this paper
is $\M$. It is a non-projective, smooth quasi-projective variety of 
dimension $6g-6$. It was first introduced in \cite{hitchin1} and 
then an algebro-geometric approach was given in \cite{nitsure} and in
\cite{simpson}.

Similarly to (\ref{split}) we have a $G$-action on $\tM$ and on
$T^*_\J$ such that:
$$\tM=(\M\times T^*_\J)/G.$$ This on the level of cohomology gives 
\begin{eqnarray}
H^*(\tM)\cong (H^*(\M))^G\otimes H^*(T^*_\J)\cong (H^*(\M))^G
\otimes H^*(\J) .\label{splitM}
\end{eqnarray}
In the case of $\M$ however we do not have the triviality of the action 
of $G$ on $H^*(\M)$, but nevertheless 
the cohomology ring of $\tM$ is determined
by the ring $(H^*(\M))^G$. 

There is quite little known about the
ring $H^*(\M)$. 
The Poincar\'e polynomial of it is calculated in \cite{hitchin1}. From
that calculation we can easily calculate the Poincar\'e polynomial of
$(H^*(\M))^G$. Nothing else is known about $H^*(\M)$. We will
return to this problem in a forthcoming paper \cite{hausel-thaddeus}. 

In this paper we settle another problem concerning the topology 
of $\M$. We calculate all intersection numbers of $\M$. Because $\M$
is non-compact we have to work with compactly supported cohomology. 
Moreover there is no compactly supported cohomology below the middle 
dimension $6g-6$. Thus the only interesting intersection numbers come
from the intersection form on $H^{6g-6}_{cpt}(\M)$. This space is $g$ 
dimensional and generated
by the compactly supported cohomology classes of the components of the
nilpotent cone, which is the zero fibre of the Hitchin map or in other
words the locus of stable Higgs bundles with nilpotent Higgs field 
(cf. Corollary 5.4 of \cite{hausel}). By considering the 
virtual 
Dirac bundle which is the analogue of the virtual 
Mumford bundle we will prove in the
last section of this paper that the ordinary cohomology classes of the
components of the nilpotent cone are trivial. 
This shows that the intersection form on $H^{6g-6}_{cpt}(\M)$ is trivial,
which is equivalent to Theorem~\ref{main}. 

As a conclusion 
it can be said that the analogue
of the Mumford relations for the moduli space of Higgs bundles 
is Theorem~\ref{main}.

\section{Hypercohomology}
\label{hyper}

In this section we recall the notion of hypercohomology of a complex from 
\cite{griffiths-harris}, and list some properties of it, which we will
use later. 

\begin{definition} Let $${\cal A}=(A_0\stackrel{d}{\longrightarrow}A_1
\stackrel{d}{\longrightarrow}A_2\longrightarrow...)$$
be a complex of coherent sheaves $A_i$ over an algebraic variety $X$. For a
covering $\underline{U}=\{U_\alpha\}$ of $X$ and each
$A_i$ we get the \v Cech cochain complex with boundary operator $\delta$:
$$(C^0(\underline{U},A_i)\stackrel{\delta}{\longrightarrow} 
C^1(\underline{U},A_i)\stackrel{\delta}{\longrightarrow}...).$$
Clearly $d$ induces operators 
$$(C^j(\underline{U},A_i)\stackrel{d}{\longrightarrow} 
C^j(\underline{U},A_j)),$$
satisfying $\delta^2=d^2=d\delta+\delta d=0$: and hence gives rise
to a double complex $$\{C^{p,q}=C^p(\underline{U},A_q);\delta,d\}.$$
The hypercohomology of the complex $\cal A$ is given by the cohomology
of the total complex of the double complex $C^{p,q}$:
$$\Hy^*(X,{\cal A})=\lim_{\underline{U}}H^*(C^*(\underline{U}),D).$$

Moreover if $\cal A$ is a complex over $X$ and $f:X\rightarrow Y$ is 
a projective morphism then for every non-negative integer $i$ 
define the sheaf $\R^i f_*({\cal A})$ over $Y$
by $$\R^i f_*({\cal A})(U)=\Hy^i(f^{-1}(U),{\cal A}).$$ 

Finally, define the pushforward of a complex to be:

$$f_!({\cal A})=
\R^0 f_*({\cal A})-\R^1 f_*({\cal A})+\R^2 f_*({\cal A})-\dots\in K(Y).$$
\end{definition}

\begin{remark} 
In this paper we will work only with two term complexes. 

There is one important property of hypercohomology which we will
make constant use of. If $$0\rightarrow {\cal A}\rightarrow
{\cal B} \rightarrow {\cal C}\rightarrow 0$$ 
is a short exact sequence of complexes
then there is a long exact sequence of hypercohomology vector spaces:
\begin{eqnarray} 0\rightarrow \Hy^0(X,{\cal A})
\rightarrow \Hy^0(X,{\cal B})\rightarrow \Hy^0(X,{\cal C})
\rightarrow \Hy^1(X,{\cal A}) \rightarrow 
\dots 
\end{eqnarray}
As an example consider the short exact sequence of two term
complexes: 
$$
\begin{array}{c}0\longrightarrow 0\\
		\uparrow \hskip1cm \uparrow \\
		 0\stackrel{0}{\longrightarrow} A_2 \\
		 \hskip.36cm \uparrow  \hskip1cm \uparrow 
			{\mbox{\scriptsize{$\cong$}}}\\
		A_1\stackrel{d}{\longrightarrow} A_2\\
		\mbox{\scriptsize{$\cong$}}\uparrow\hskip1cm\uparrow 
					\hskip.36cm \\
		A_1\stackrel{0}{\longrightarrow} 0\\
		\uparrow\hskip1cm\uparrow\\
		0\longrightarrow 0
\end{array}		
$$
The long exact sequence in this case is:
\begin{eqnarray}0\rightarrow \Hy^0(X,{\cal A})
\rightarrow H^0(X,A_1)\rightarrow H^0(X,A_2)\rightarrow 
\Hy^1(X,{\cal A}) \rightarrow \dots \label{longexact}
\end{eqnarray}
which we will call the hypercohomology long exact sequence of the two
term complex ${\cal A}=A_1\stackrel{d}\rightarrow A_2$.

Consequently if ${\cal A}=A_1\stackrel{d}{\rightarrow} A_2$
is a two term complex over $X$ and $f:X\rightarrow Y$ is a projective 
morphism
then we have:

\begin{eqnarray} 0\rightarrow  \R^0f_*(X,{\cal A})
\rightarrow R^0f_*(X,A_1)\rightarrow R^0f_*(X,A_2)\rightarrow 
\R^1f_*(X,{\cal A}) \rightarrow
\dots, \label{long} \end{eqnarray}
a long exact sequence of sheaves over $Y$.
\end{remark}

\section{A vanishing theorem}

\label{vanishin}

\begin{definition} 
The complex $\compE$ with $E$ a
vector bundle on $\Sigma$, $K$ the canonical bundle of $\Sigma$, 
and $\Phi\in H^0(\Sigma,\Hom(E,E\otimes K))$, 
is called a {\em Higgs bundle}, while $\Phi$ is called the 
{\em Higgs field}. 

We define a morphism $\Psi:{\cal E}_1 \rightarrow {\cal E}_2$
between two Higgs bundles \newline
\mbox{${\cal E}_1=E_1\stackrel{\Phi_1}{\rightarrow} E_1\otimes K$} and  
\mbox{${\cal E}_2=E_2\stackrel{\Phi_2}{\rightarrow} E_2\otimes K$} 
to be a homomorphism of vector bundles 
$\Psi\in \Hom(E_1,E_2)$ such that the following diagram commutes:
$$
\begin{array}{c}
 \  E_1\stackrel{\Phi_1}{\longrightarrow} E_1\otimes K \ \\
\hskip.35cm\mbox{\scriptsize{$\Psi$}}\downarrow\hskip1.8cm \downarrow 
\mbox{\scriptsize{$\Psi\otimes id_K$}} \\
 \  E_2\stackrel{\Phi_2}{\longrightarrow} E_2\otimes K\ 
\end{array}
$$ 

Moreover we say that 
${\cal E}_1$ is a {\em Higgs
subbundle} of ${\cal E}_2$ if 
$\Psi\in \Hom(E_1,E_2)$ is injective and a morphism of Higgs
bundles. We denote this by ${\cal E}_1\subset {\cal E}_2$. 
In this case we can easily construct  
the quotient Higgs bundle ${\cal E}_2/{\cal E}_1$ together with a
surjective morphism of Higgs bundles $\pi:{\cal E}_2\rightarrow {\cal
E}_2/{\cal E}_1$ whose kernel is exactly ${\cal E}_1$.
\label{higgs}
\end{definition}

\begin{remark} It is a tautology that morphisms of Higgs bundles
form the hypercohomology\footnote{In connection
with Higgs bundles the language of 
hypercohomology was first used in \cite{simpson}.
In \cite{biswas-ramanan} 
it was used to describe the tangent space to $\M$.} 
vector space  
$\Hy^0(\Sigma, E_1^*\otimes E_2\stackrel{\left[\Phi_1,\Phi_2\right]}
{\longrightarrow}
E_1^*\otimes E_2\otimes K)$ where the homomorphism 
$\left[\Phi_1,\Phi_2\right]$ is given
by: ${\left[\Phi_1,\Phi_2\right](\Psi):=(\Psi\otimes
id_K)\Phi_1-\Phi_2\Psi}$ for $\Psi\in \Hom(E_1,E_2)$. 
\end{remark}

Now we can define the notion of stability of Higgs bundles:

\begin{definition} If $E$ is a vector bundle over $\Sigma$ then its
{\em slope} is defined by $\mu(E):=\deg(E)/\rank(E)$. The slope 
$\mu(\cal E)$ of a Higgs
bundle\newline $\cE=\compE$ 
is defined as the slope $\mu(E)$ of its vector bundle $E$. Now a Higgs
bundle is called {\em stable} if it has strictly larger slope than any of its
proper Higgs subbundles. 
\end{definition} 

The main result of this section is the following theorem, the second
part of which is Proposition (3.15) in \cite{hitchin1}:

\begin{theorem} Let ${\cal E}=E\stackrel{\Phi}{\rightarrow}E\otimes K$
and ${\cal F}=F\stackrel{\Psi}{\rightarrow}{F\otimes K}$ be stable Higgs
bundles with $\mu({\cal F})<\mu({\cal E})$. Then the only morphism
from ${\cal E}$ to ${\cal F}$ is the trivial one. In other words 
$$\Hy^0(\Sigma, E^*\otimes F\stackrel{[\Phi,\Psi]}{\longrightarrow}
E^*\otimes F\otimes K)=0.$$

Moreover if $\mu({\cal F})=\mu({\cal E})$, then there is a non-trivial
morphism $f:{\cal E} \rightarrow {\cal F}$ if and only if 
${\cal E}\cong {\cal F}$ in which case every non-trivial morphism 
$ f:{\cal E} \rightarrow {\cal F}$ is an isomorphism and 
\begin{eqnarray}
\dim(\Hy^0(\Sigma, E^*\otimes F\stackrel{[\Phi,\Psi]}{\longrightarrow}
E^*\otimes F\otimes K))=1.
\label{iso}
\end{eqnarray}
\label{van}
\end{theorem}

\begin{proof} For the proof we need a lemma of Narasimhan and Seshadri
(cf. section 4 in \cite{narasimhan-seshadri}):
\begin{lemma}
Let $E$ and $F$ be two vector bundles over the Riemann surface
$\Sigma$ with a non-zero homomorphism $f:E\rightarrow F$, then $f$ has
the following canonical factorisation: 
$$
\begin{array}{c}
 0\longrightarrow E_1 \longrightarrow E
 \stackrel{\eta}{\longrightarrow} E_2 \longrightarrow 0 \\
\hskip1.65cm \downarrow {\mbox{\scriptsize{f}}} \hskip.75cm 
\downarrow {\mbox{\scriptsize{g}}} \\
0 \longleftarrow F_2 \longleftarrow F \stackrel{i}{\longleftarrow} F_1
 \longleftarrow 0
\end{array}
$$
where $E_1,E_2,F_1$ and $F_2$ are vector bundles, each row is exact, 
$f=ig\eta$ and $g$ is of maximal rank, i.e. $\rank(E_2)=\rank(F_1)=n$ and 
$\Lambda^n(g):\Lambda^n(E_2)\rightarrow \Lambda^n(F_1)$ is a non-zero 
homomorphism. In other words $g$ is an isomorphism on a Zariski open
subset
$U$ of $\Sigma$. $F_1$ is called the subbundle of $F$ generated by the
image of $f$. $\square$
\label{narasimhan-seshadri}
\end{lemma}

Let $f:{\cal E}\rightarrow {\cal F}$ be a non-zero morphism of Higgs
bundles. In particular $f:E\rightarrow F$ is a homomorphism of
vector bundles. 

Construct the canonical factorisation of $f$ of the
above lemma. Consider the Zariski open subset $U$ of $\Sigma$ where
$g$ is an isomorphism. Here clearly $\ker(f\mid_U)=\ker(\eta\mid_U)
= E_1\mid_U$. Now $\ker(f\mid_U)$ being the kernel of a morphism of
Higgs bundles is $\Phi$-invariant, i.e. a Higgs subbundle of ${\cal E}\mid_U$. 
Thus $E_1\mid_U$ is a Higgs subbundle of ${\cal E}\mid_U$. This means
that $\Phi(E_1)$ is contained in $E_1\otimes K\subset E\otimes K$ on
$U$. Because $U$ is Zariski open in $\Sigma$ it follows that 
${\cal E}_1=E_1\stackrel{\Phi}{\rightarrow} E_1\otimes K$ is
a Higgs subbundle of ${\cal E}$. Let ${\cal E}_2=E_2
\stackrel{\tilde{\Phi}}{\rightarrow} E_2\otimes K$ denote the quotient
Higgs bundle.  

Similarly $\im(\alpha)\mid_U= F_1\mid_U$ is $\Psi$-invariant, thus
${\cal F}_1=F_1\stackrel{\Psi}{\rightarrow} F_1\otimes K$ is a
Higgs subbundle of $F$. 

By assumption $\mu({\cal F)}<\mu({\cal E})$, 
stability of ${\cal E}$ gives $\mu({\cal E})\leq \mu ({\cal E}_2)$ (it
may happen that $E=E_2$) and
because $g$ is of maximal rank we get 
$\mu({\cal E}_2)=\mu(E_2)\leq \mu(F_1)=\mu({\cal F}_1)$.
Thus $\mu({\cal F})< \mu({\cal F}_1)$ contradicting the stability of
$F$. 

If $\mu(\cE)=\mu({\cal F})$ then the above argument leaves the only
possibility that $\eta$, $g$ and $i$ are isomorphisms, showing that
$f$ must be an isomorphism. Suppose that we have such an isomorphism
$f$ of Higgs bundles. Then consider $h:\cE \rightarrow {\cal F}$ 
another non-zero morphism of Higgs bundles. In particular $h\in 
\Hom(E,F)$. Let $\lambda$ be an eigenvalue of the homomorphism 
$f_p^{-1}h_p\in\Hom(E_p,E_p)$. Then the homomorphism $h-\lambda f$ is
not an isomorphism, though clearly a morphism of Higgs bundles. From
the above argument this means that $h-\lambda f=0$.     

The result follows. 
\end{proof}

\begin{corollary} 
For any stable Higgs bundle ${\cal E}$ with $\mu({\cal E})<0$:
 \begin{eqnarray}\Hy^0(\Sigma,{\cal E})=0,\label{h0}\end{eqnarray}
for any stable Higgs bundle ${\cal E}$ with $\mu({\cal E})>0$:
\begin{eqnarray} \Hy^2(\Sigma,{\cal E})=0.\label{h2} \end{eqnarray}
If $\cE$ is a stable Higgs bundle with $\mu(\cE)=0$ and $\cE \ncong
\cE_0=\calO_{\Sigma}\stackrel{0}{\rightarrow}\calO_\Sigma \otimes K$
then both (\ref{h0}) and (\ref{h2}) hold. 
\label{vanish}
\end{corollary}
\begin{proof} 
For the first part consider the Higgs bundle 
${\cal E}_0=\calO_{\Sigma}\stackrel{0}{\rightarrow}\calO_\Sigma \otimes
K$. 
Being of rank $1$ it is obviously stable, with $\mu({\cal
E}_0)=0$. Now the previous theorem yields that there are no nontrivial
morphisms
from ${\cal E}_0$ to ${\cal E}$, which in the language of
hypercohomology is exactly $\Hy^0(\Sigma,{\cal E})=0$, which we had to
prove. 

For the second part Serre duality gives that 
$\Hy^2(\Sigma,\cE)\cong (\Hy^0(\Sigma,\cE^*\otimes K))^*$. Now 
 clearly $\cE^*\otimes K$ is stable and 
$\mu(\cE^*\otimes K)=-\mu(\cE)<0$. Thus the first part gives the
second. 

Likewise, the third statement follows by referring to the last part of 
Theorem~\ref{van}.
\end{proof}

\section{Universal bundles}
\label{univbundle}

Nitsure showed that $\tM$ is a coarse moduli space. 
Here we show that 
$\tM$ is in fact a {\em fine} moduli space. We closely follow
the proof of Theorem 5.12 in \cite{newstead} and (1.19) of
\cite{thaddeus3}. All the ingredients have already appeared in the 
unpublished \cite{thaddeus1}. 

\begin{definition} Two families $\cE_T$ and $\cE_T^\prime$ 
of stable Higgs bundles over $T\times \Sigma$ are said to be 
{\em equivalent}, (in symbols $\cE_T\sim \cE_T^\prime$) 
if there exists a line
bundle $L$ on $T$ such that 
 $\cE^\prime_T\cong\cE_T\otimes\pi_T^*(L)$. 
\end{definition}

The next lemma, which is taken from \cite{thaddeus1}, 
shows that two families are equivalent iff they give rise to the same
map to the coarse moduli space $\tM$.

\begin{lemma} If 
$\cE_T=\E_T\stackrel{\uPhi}{\rightarrow}\E_T\otimes
\K$ and $\cE_T^\prime=
\E^\prime_{T}\stackrel{\uPhi^\prime}{\rightarrow}
\E^\prime_T\otimes
\K$ are families of stable
Higgs bundles
over $T\times \Sigma$ such 
that \begin{eqnarray}\cE_T\mid_{\{t\}\times \Sigma}\cong
\cE^\prime_T\mid_{\{t\}\times \Sigma}\label{con}
\end{eqnarray} for each $t\in T$, then $\cE_T\sim \cE_T^\prime$.
\label{ramanan}
\end{lemma}
  
\begin{proof} Let ${\cal F}:= \E^*_T\otimes\E_T^\prime
\stackrel{[\uPhi_T,\uPhi_T^\prime]}{\longrightarrow} 
\E^*_T\otimes\E_T^\prime\otimes\K$.      
We define $L=\R^0{\pi_T}_*({\cal F})$. By (\ref{con}) and 
(\ref{iso}) this is a line
bundle over $T$. By the projection formula 
the sheaf $\R^0{\pi_T}_*({\cal F}\otimes 
\pi_T^*(L^*))$ is just ${\calO}_T$, the structure sheaf. A
non-zero section ${\mathbf \Psi}\in H^0(T,\R^0{\pi_T}_*({\cal F}\otimes 
\pi_T^*(L^*)))$ for every $t\in T$ gives ${\mathbf
\Psi}\mid_{\{t\}\times \Sigma}:
(\cE_T\otimes \pi_T^*(L))\mid_{\{t\}\times \Sigma} \rightarrow 
\cE^\prime_T\mid_{\{t\}\times \Sigma}$ a non-zero morphism of  
Higgs bundles, which is by Theorem~\ref{van} an isomorphism. 

The result follows.
\end{proof}

Now we prove the existence of universal Higgs bundles (cf. \cite{thaddeus1}):

\begin{proposition} Universal Higgs bundles $\cE_\tM=\ucompE$ over
$\tM\times \Sigma$ do exist. 
\label{universal}
\end{proposition}

\begin{proof} The proof is analogous to the proof of Theorem 5.12 of
\cite{newstead} using the GIT construction of Nitsure \cite{nitsure}
(cf. also (1.19) of \cite{thaddeus3}).  

First we recall the construction of $\tM_{2k-1}$ from \cite{nitsure}. 
Let $n=2k-1+2(1-g)$ with
$k$ large enough. Then by Corollary 3.4 of \cite{nitsure} for any stable
Higgs bundle $\compE$, $E$ is a quotient of $\calO^n_\Sigma$. Let $\Q$
be the quot scheme of all quotient sheaves $\calO^N_\Sigma\rightarrow
{\cal F}$ of rank $2$ and degree $d$. Let $\calO^n_{\Sigma\times
\Q}\rightarrow U$ be the universal quotient sheaf on $\Sigma\times
\Q$. Let $R\subset \Q$ be the subset of all $q$ for which ${\cal F}_q$
is locally free and the map $H^0(\Sigma,\calO_\Sigma^n)\rightarrow
H^0(\Sigma,U_q)$ is an isomorphism. 

It follows from Proposition 3.6 of \cite{nitsure} that there
exists a locally universal family for stable Higgs bundles of degree
$2$ and degree $2k-1$ given by
$\cE_s=\E_s\stackrel{\uPhi_s}{\longrightarrow}\E_s \otimes \K$ over 
$F_s\times \Sigma$ where $F_s$ is an open subset of a linear
$R$-scheme 
$F\rightarrow R$ and
$\cE_s=\cE_F\mid_{F_s\times \Sigma}$ where
$\cE_F=\E_F\stackrel{\uPhi_F}{\longrightarrow} \E_F\otimes \K$ is a
family of Higgs bundles over $F$.

First by Theorem 5.3 of \cite{newstead} $GL(n)$ acts on $R$. 
Now $GL(n)$ acts equivariantly on the $R$-scheme
$F\rightarrow R$, which gives a $GL(n)$ equivariant complex $\cE_F$. 
The centre of
$GL(n)$ acts trivially on $F$ and by multiplication on
$\cE_F$. Nitsure constructs $\tM_{2k-1}$ in Theorem 5.10 of 
\cite{nitsure} as a good quotient of
$F_s$ by $PGL(n)\cong GL(n)/Z(GL(n))$. 

The proof of Lemma 5.11 of \cite{newstead} gives a
$GL(n)$-equivariant line bundle $L$ over $R$ (although in Lemma 5.11
of \cite{newstead} $L$ is constructed only over $R_s$ the same
construction works over the whole $R$) for which $Z(GL(n))$ acts
on $L$ by scalar multiplication. Now for the $GL(n)$-equivariant bundle
$\cE_F\otimes (\pi_F\circ g)^*(L^{-1})$ the centre acts trivially thus
it descends to a $PGL(n)$-equivariant complex over $F\times \Sigma$. 
This gives a $PGL(n)$-equivariant locally universal family  
$\cE_s\otimes (\pi_F\circ g)^*(L^{-1})$ 
over
$F_s\times \Sigma$. By Kempf's descent lemma (cf. Theorem 2.3 of 
\cite{drezet-narasimhan})) the $PGL(n)$-equivariant bundle 
$\E_s\times (\pi_F\circ g)^*(L^{-1})$ descends to a bundle to 
the good quotient $\tM_{2k-1}\times \Sigma$ and since  the section 
$\uPhi_s$ is invariant, it also descends.     
Clearly the resulting complex $\cE_{\tM_{2k-1}}$ 
then will be a universal Higgs bundle over $\tM_{2k-1}$. 
(A similar situation appears in (1.19) of \cite{thaddeus3}.)

Finally from a universal Higgs bundle over $\tM_{2k-1}$ one can easily
construct universal Higgs bundles over any $\tM_{2l-1}$. 

The result follows.
\end{proof}

As in Theorem 5.12 of \cite{newstead} and (1.19) of \cite{thaddeus3} 
our Lemma~\ref{ramanan} and Proposition~\ref{universal} gives:

\begin{corollary} The space $\tM$ is a fine moduli space for 
rank $2$ stable Higgs bundles of degree $1$ with respect to the
equivalence $\sim$ of families of stable Higgs bundles.   
\end{corollary}

As another  consequence of Proposition~\ref{universal} and 
Lemma~\ref{ramanan} we see that although $\E_\tM$ is not unique
$\End(\E_\tM)$ is.  Moreover it is clear that 
by setting $\E_\M=\E_\tM\mid_{\M\times \Sigma}$ we have
\begin{eqnarray}
c(\End(\E_\tM))=c(\End(\E_\M))\otimes 1\label{c}
\end{eqnarray}
in the decomposition (\ref{splitM}).
 
Thus from the K\"unneth decomposition of $\End(\E_\M)$ we get 
universal classes
$$c_2(\End(\E_\M))=2\alpha_\M\otimes \sigma+\sum_{i=1}^{2g}4\psi^i_\M 
\otimes e_i -
\beta_\M\otimes 1$$
in $H^4(\M\times \Sigma)\cong
\sum_{r=0}^{4} H^r(\M)\otimes H^{4-r}(\Sigma)$
for some $\alpha_\M\in H^2(\M)$, $\psi^i_\M\in H^3(\M)$ and $\beta_\M\in
H^4(\M)$. 

Clearly $\alpha_\M\mid_\N=\alpha_\N$, $\psi^i_\M\mid_\N=\psi^i_\N$ and
$\beta_\M\mid_\N=\beta_\N$. 

Though $\E_\M$ is not unique we can still write its Chern classes 
in the K\"unneth decomposition 
(cf. proof of Newstead's theorem in \cite{thaddeus2}), 
getting $c_1(\E_\M)=1\otimes \sigma +\beta_1\otimes 1$,
where
$\beta_1\in H^2(\M)$
(note that $\M$ being simply connected by \cite{hitchin1} $H^1(\M)=0$)
and $c_2(\E_\M)=\alpha_2 \otimes \sigma + \sum_{i=1}^{2g}a_i \otimes
e_i + 
\beta_2\otimes 1$,
where $\alpha_2\in H^2(\M)$, $a_i\in H^3(\M)$ and $\beta_2\in
H^4(\M)$. 
Because
$4c_2(\E_\M)-c_1^2(\E_\M)=c_2(\End(\E_\M))$, we get
$\alpha_\M=2\alpha_2-\beta_1$ and $\beta=\beta_1^2-4\beta_2$.
Because $\Pic(\M)\cong H^2(\M,\Z)$ (cf. \cite{hausel}) we can
normalize $\E_\M$ uniquely such that $\beta_1=\alpha_\M$ i.e. 
$c_1(\E_\M)=1\otimes \sigma +\alpha_\M\otimes 1$.    

\begin{definition} The universal Higgs bundle $\cE_\M$ is
normalized if\newline $c_1((\E_\M)_p)=\alpha_\M$, where 
$(\E_\M)_p:=\E_\M\mid_{\M\times \{p\}}$.
\end{definition}

We also need to work out the Chern classes of $\E_\tM$.
It is easy to
see that $c(\E_\tM)$ in the decomposition (\ref{splitM}) is the product
of $c(\E_\tM)\mid_{\M\times \Sigma}$ and $c({\bL}_1)$, where 
${\bL}_1$ is some universal line bundle over $\J\times \Sigma$. 

\begin{definition} We call the universal
Higgs bundle $\cE_\tM$ {\em normalized} if in the decomposition 
(\ref{splitM}) 
\begin{eqnarray} 
c_1((\E_\tM)_p)=\alpha_\M \label{c1},
\end{eqnarray}
where $(\E_\tM)_p=\E_\tM\mid_{\tM\times \{p\}}.$
\end{definition}

\begin{remark} Since
$4c_2((\E_\tM)_p)-c_1((\E_\tM)_p)^2=c_2\left(\End((\E_\tM)_p)\right)$,
for a normalized universal Higgs bundle over $\tM\times \Sigma$
(\ref{c}) and (\ref{c1}) yield:
\begin{eqnarray}c_2((\E_\tM)_p)=\frac{(\alpha_\M^2-\beta_\M)}{4}
  \label{c2}
\end{eqnarray} 
\end{remark}

Finally, given a universal Higgs bundle $\cE_\tM$ over $\tM\times
\Sigma$, we introduce
a {\em universal Higgs bundle of degree $2k-1$} by setting 
$$\cE_{\tM}^k:=\cE_{\tM}\otimes \pi_\Sigma^*(L_p^{k-1}),$$ where $L_p$ is
the line bundle of the divisor of the point 
$p\in \Sigma$. It is called normalized
if $\cE_\tM$ is normalized. As a matter of fact 
$\cE_{\tM}^k$ can be thought of as a pull back of 
a universal Higgs bundle from $\tM_{2k-1}\times \Sigma$.

\section{The virtual Dirac bundle, $\D_k$}
\label{virdirac}

The strategy of the proof of Theorem~\ref{main} will be to examine the
virtual Dirac bundle  $\D_k$ which is defined in the following:

\begin{definition} The {\em virtual Dirac bundle} is\footnote{Recall
the
definition of the pushforward of a complex from Section~\ref{hyper}.} 
$$\D_k:=-{\pi_\tM}_!(\cE_{\tM}^k)\in K(\tM),$$ 
where $\cE^k_\tM$ is a normalized
universal Higgs bundle of degree $2k-1$ 
and $\pi_\tM:\tM\times \Sigma\rightarrow \tM$ is
the projection to $\tM$.
\end{definition}

The name is justified by Hitchin's construction
\cite{hitchin2}\footnote{Cf. Subsection 1.1.5 of \cite{hausel2}.} of
$\D_k$ related to the space of solutions of an equation on $\Sigma$, 
which is locally 
the dimensional reduction of the Dirac equation in $\R^4$ coupled to
a self-dual Yang-Mills field. 

The
virtual Dirac bundle is  
a priori $${\pi_\M}_!(\cE_{\tM}^k)=
-\R^0{\pi_\tM}_*(\cE_{\tilde{\M}}^k)+
\R^1{\pi_\tM}_*(\cE_{\tilde{\M}}^k)-
\R^2{\pi_\tM}_*(\cE_{\tilde{\M}}^k)$$ a
formal sum of three coherent sheaves. Corollary~\ref{vanish} shows
that one of these sheaves always vanishes: 
if $k> 0$, then $\R^2=0$, if $k\leq 0$ then $\R^0=0$. From now on $k$ is
assumed to be positive.

In this section we show that we
can think of the virtual Dirac bundle as the virtual 
degeneracy sheaf of a homomorphism
of vector bundles. More precisely we prove:

\begin{theorem}
There exist two vector bundles $V$ and $W$ over $\tM$ together with a 
homomorphism $f:V\rightarrow W$ of vector bundles, 
whose kernel and cokernel are respectively  
$\R^0{\pi_\M}_*(\cE_{\tilde{\M}}^k)$ and 
$\R^1{\pi_\M}_*(\cE_{\tilde{\M}}^k)$. In other
words there is an exact sequence of sheaves:
$$0\rightarrow \R^0{\pi_\tM}_*(\cE_{\tM}^k)
\rightarrow V\stackrel{f}{\rightarrow} W \rightarrow
\R^1{\pi_\tM}_*(\cE_{\tM}^k)\rightarrow
0.$$
\label{deg}
\end{theorem}
\noindent {\it Proof\footnote{The 
idea of the proof was suggested by Manfred Lehn.}.}
First we need a lemma. 

\begin{lemma} Let $X$ be a smooth quasi-projective variety 
and $\Sigma$ a
smooth projective curve. 
If $E$ is a locally free sheaf over $X\times {\Sigma}$ 
then there exists a vector bundle $F$ over $X \times \Sigma$ 
with a surjective vector
bundle homomorphism $g_E:F\rightarrow E$ such that 
$R^0{\pi_{X}}_*(F)=0$. 
We will call $F$ a {\em sectionless resolution} of $E$.
\end{lemma}

\begin{proof} The lemma is a special case of Proposition 2.1.10 of 
\cite{huybrechts-lehn}. We have to only note that $X$ as an
algebraic variety 
is a $\C$-scheme of finite type, ${\pi_{X}}_*:{X}\times
\Sigma\rightarrow X$ is clearly a smooth projective morphism of 
relative dimension
$1$ and $E$ being locally free is flat over $X$. 

The proof is rather simple so we sketch it here.
Let us denote by $E_x$ the vector bundle $E\mid_{\{x\}\times \Sigma}$
over $\Sigma$. Fix an ample line bundle $L$ on $\Sigma$. Then it is
well known that for big enough $k$ the vector bundle $E_x\otimes L^k$ 
is generated by its sections and $H^1(\Sigma;E_x\otimes L^k)=0$. 
Let us denote by $X_k\subset X$ those
points $x$ for which $E_x\otimes L^k$ is generated by its sections and
$H^1(\Sigma;E_x\otimes L^k)=0$. It
is standard that $X_k$ is a Zariski open subset of $X$. Thus we have a
covering $X=\bigcup X_k$ of $X$ by Zariski open subsets. It is again
standard that the Zariski topology of an algebraic variety is
noetherian\footnote{Cf. Example 3.2.1 on p. 84 of \cite{hartshorne}.},
which yields that we have some $k$ such that $X_k=X$. It is
now immediate that $$F=\pi_{\Sigma}^*(L^{-k})\otimes
\pi_X^*\left( (\pi_X)_*(E\otimes \pi_\Sigma^*(L^k))\right)$$ has the required
properties.

The result follows.
\end{proof} 

\begin{proposition} Let ${\Sigma}$ be a smooth projective curve 
and $X$ be a smooth quasi-projective variety. Let 
$\cE=E\stackrel{f}{\rightarrow} F$ be a complex of vector bundles
on $X\times {\Sigma}$. Let $g_F:A\rightarrow F$ be 
a sectionless resolution of $F$. Let $M$ be the fibred product of 
$f$ and $g_F$. This comes with projection maps $p_F:M\rightarrow F$ and
$p_{A}:M\rightarrow A$. Let $g_M:A_2\rightarrow M$ 
be a sectionless resolution of $M$, and denote $j=g_M \circ p_{A_2}$. Finally, let
$A_1=\ker{g_M}$ and $i:A_1\rightarrow A_2$ the embedding. 
The situation is shown 
in the following diagram:
$$
\begin{array}{cc}
&E\hskip.2cm\stackrel{f}{\longrightarrow}\hskip.2cm F\\
&\hskip.3cm\nwarrow\hskip1cm\\
&\ \hskip.5cm M \hskip.4cm\uparrow\\
&\nearrow\hskip.4cm\searrow\\
0\longrightarrow A_1\stackrel{i}{\longrightarrow}& A_2\hskip.2cm
\stackrel{j}{\longrightarrow} \hskip.2cm A
\end{array}
$$
In this case 
the cohomology of the complex 
$$R^1{\pi_X}_*(A_1)\stackrel{i_*}{\longrightarrow} 
R^1{\pi_X}_*(A_2) \stackrel{j_*}{\longrightarrow} 
R^1{\pi_X}_*(A)$$
calculates the sheaves $\R^0{\pi_X}_*(\cE)$,
 $\R^1{\pi_X}_*(\cE)$ and $\R^2{\pi_X}_*(\cE)$
respectively. 
In other
words 
\begin{eqnarray}\R^0{\pi_X}_*(\cE)&\cong& \ker(i_*)\label{r0}\\  
\R^1{\pi_X}_*(\cE)&\cong &\ker(j_*)/\im(i_*) \label{r1}\\
\R^2{\pi_X}_*(\cE)&\cong& \coker(j_*)\label{r2}.
\end{eqnarray}
\label{diagram}
\end{proposition}
\begin{proof} Let us recall the definition of the fibred product: 
$M:=\ker(f-g_F:E\oplus A\rightarrow F)$. This comes equipped with
 two obvious projections $p_E:M\rightarrow E$ and
$p_{A}:M\rightarrow A$. Because $g_F$ is surjective
$f-g_F$ is also surjective. Thus $M$ is a vector bundle. By construction
the kernel
of $p_E$ is isomorphic to the kernel of $g_F$. Denote it by $B$. 
This says that the
following diagram is commutative and has two exact columns: 
$$
\begin{array}{c}0\longrightarrow 0\\
		\uparrow \hskip1cm \uparrow \\
		 E\stackrel{f}{\longrightarrow} F \\
		\mbox{\scriptsize{$p_E$}} \uparrow  \hskip1cm \uparrow
		\mbox{\scriptsize{$g_F$}}\\
		M\stackrel{p_{A}}{\longrightarrow} A\\
		\uparrow\hskip1cm\uparrow\\
		B\stackrel{\cong}{\longrightarrow} B\\
		\uparrow\hskip1cm\uparrow\\
		0\longrightarrow 0
\end{array}		
$$

If ${\cal A}$ denotes
the complex ${\cal A}=M\stackrel{p_{A}}{\rightarrow}A$ and ${\cal
B}$ the complex ${\cal B}=B\stackrel{\cong}{\rightarrow} B$, then the
above diagram is just a short exact sequence of complexes 
$$0\longrightarrow {\cal B}\longrightarrow {\cal A}\longrightarrow
{\cE}\longrightarrow 0.$$

Clearly $\R^i{\pi_X}_*({\cal B})$ vanishes for all $i$. 
(Any hypercohomology of an isomorphism is $0$.) Thus the long exact
sequence of the above short exact sequence gives the isomorphisms 
\begin{eqnarray}
\R^0{\pi_X}_*({\cal E})&\cong & \R^0{\pi_X}_*({\cal A})
\label{1r0}\\
\R^1{\pi_X}_*({\cal E})&\cong & \R^1{\pi_X}_*({\cal A})
\label{1r1}\\
\R^2{\pi_X}_*({\cal E})&\cong & \R^2{\pi_X}_*({\cal A})
\label{1r2}
\end{eqnarray}
 
Because $A$ is a sectionless resolution of $M$, we have 
$R^0{\pi_X}_*(A)=0$
thus the long exact sequence of the push forward of the complex $\cal
A$ breaks up into two exact sequences:
$$0\rightarrow \R^0{\pi_X}_*({\cal A})\rightarrow 
R^0{\pi_X}_*(M)\rightarrow 0,$$
and 
$$
0 \longrightarrow 
\R^1{\pi_X}_*({\cal A})\longrightarrow R^1{\pi_{\cal
W}}_*(M)\stackrel{{p_{A}}_*}\longrightarrow R^1{\pi_X}_*(A)
\longrightarrow \R^2{\pi_X}_*({\cal A})\longrightarrow 0.$$
Thus 
\begin{eqnarray}
\R^0{\pi_X}_*({\cal A})&\cong& R^0{\pi_X}_*(M)
\label{2r0}\\ 
\R^1{\pi_X}_*({\cal A})&\cong& \ker({p_{A}}_*)
\label{2r1}\\ 
\R^2{\pi_X}_*({\cal A})&\cong& \coker({p_{A}}_*). \label{2r2}
\end{eqnarray}

Now consider the short exact sequence:
$$0\longrightarrow A_1\stackrel{i}{\longrightarrow} A_2 
\stackrel{g_M}{\longrightarrow} M\longrightarrow 0. \label{exact1}$$
$R^0{\pi_X}_*(A_2)=0$ because $A_2$ is a 
sectionless resolution of $M$ and hence we get the exact
sequence of sheaves:
\begin{eqnarray}0&\longrightarrow& R^0{\pi_X}_*(M)\longrightarrow 
R^1 {\pi_X}_*(A_1)\stackrel{i_*}{\longrightarrow} 
R^1{\pi_X}_*(A_2)
\stackrel{{g_M}_*}\nonumber\\&{\longrightarrow}&R^1{\pi_{\cal
W}}_*(M)\longrightarrow 0.
\label{exact}
\end{eqnarray}

Thus $\ker(i_*) \cong R^0{\pi_X}_*(M)$ which by (\ref{2r0}) and
(\ref{1r0}) proves (\ref{r0}). 

Since ${g_M}_*$ is a surjection $\coker(j_*)\cong
\coker({p_{A}}_*)$. This together with (\ref{2r2}) and (\ref{1r2})
gives (\ref{r2}).

Finally, consider the commutative diagram:
$$
\begin{array}{c}R^1{\pi_{\cal
W}}_*(M)\stackrel{\cong}{\longrightarrow} 
R^1{\pi_X}_*(M)\\
\mbox{\scriptsize{${g_M}_*$}}\uparrow \hskip1.5cm \downarrow
\mbox{\scriptsize{${p_{A}}_*$}}\\
R^1{\pi_X}_*(A_2)\stackrel{j_*}{\longrightarrow}R^1{\pi_X}_*(A)
\end{array}
$$
Since ${g_M}_*$ surjective by (\ref{exact}) we get that $\ker(j_*)/
\ker({g_M}_*)\cong \ker({p_{A}}_*)$. From (\ref{exact}) clearly 
$\ker({g_M}_*)\cong \im(i_*)$, thus $\ker(j_*)/\im(i_*)\cong 
\ker({p_{A}}_*)$. This together with (\ref{2r1}) and (\ref{1r1})
proves (\ref{r1}).
\end{proof}

\begin{corollary} If $\R^2{\pi_X}_*(\cE)=0$, 
in the situation of Proposition~\ref{diagram},
 then there exist two vector bundles $V$ and
$W$ over $X$ together with a homomorphism $f:V\rightarrow W$,
whose kernel and cokernel are $\R^0{\pi_X}_*(\cE)$ and 
$\R^1{\pi_X}_*(\cE)$ respectively. I.e. the following sequence
is exact:
$$0\rightarrow \R^0{\pi_X}_*(\cE)
\rightarrow V\stackrel{f}{\rightarrow} W \rightarrow
\R^1{\pi_X}_*(\cE)\rightarrow
0.$$
\label{cor}
\end{corollary}

\begin{proof} From the long exact sequence corresponding to
(\ref{exact1}), we have \newline $R^0{\pi_X}_*(A_1)=0$. Let $V$ be the
vector bundle $R^1{\pi_X}_*(A_1)$.

Moreover $R^1{\pi_X}_*(A_2)$ and $R^1{\pi_X}_*(A)$ are
also vector bundles because $A_2$ and $A$ are sectionless
resolutions. Furthermore the assumption $\R^2{\pi_X}_*(\cE)=0$
shows that $j_*$ is surjective. Let $W$ be the vector bundle
$\ker(j_*)$, and $f$ be the map $i_*:V\rightarrow W$. 

The result follows from  Proposition~\ref{diagram}.
\end{proof}

The proof of Theorem~\ref{deg} is completed by 
Corollary~\ref{cor} noting that by
Corollary~\ref{vanish} we have $\R^2{\pi_{\tM}}_*(\cE^k_{\tM})=0$.
$\square$

\section{The degeneracy locus $D_k$}
\label{degloc}

\begin{definition} The degeneracy locus $D_k:=\{\cE\in \tM: 
\Hy^0(\Sigma,\cE_\tM^k)\neq 0)\}$ 
is the locus\footnote{For a rigorous construction of degeneracy
loci cf. \cite{arbarello-et-al} p.83. Our degeneracy locus is the
$k$-th degeneracy locus of \cite{arbarello-et-al}, where
$k=\rank(V)-1$.} where $\D_k$ fails to be a
vector bundle, i.e. where $f$ of Theorem~\ref{deg} fails to be an 
injection.  
\end{definition}

The aim of this section is to give a description of the degeneracy
locus $D_k$. For this we need  a refinement  of Theorem 5.5 of \cite{hausel},
which still follows from the proof of Proposition (19) of \cite{thaddeus1}.

\begin{definition} The nilpotent cone $N\subset \M$ 
is the set of stable Higgs bundles with nilpotent Higgs field. In
other words it is $\chi^{-1}(0)$: the zero fibre of the Hitchin map. 

Similarly $\tilde{N}:=\tilde{\chi}^{-1}(0)\subset \tM$.
\end{definition}  

\begin{proposition} The nilpotent cone is a compact union of $3g-3$ 
dimensional manifolds: $$N=\N\cup\bigcup^{g-1}_{k=1}E_k,$$
where each $E_k$ is biholomorphic to the total space of a vector
bundle over $N_k$, the $k$-th component of the fixed point set of
the $\C^*$-action. 

Moreover $E_k$ can be characterised as the locus of those stable Higgs
bundles $\cE=\compE$ which have a unique 
subbundle $L_\cE$ of degree $1-k$ killed by the non-zero Higgs
field $\Phi$. 
\label{cone}
\end{proposition} 

\begin{proof}The first part is proved in Theorem 5.5 of \cite{hausel}.

For the second part consider a universal Higgs bundle $\cE_\M$
over $\M\times \Sigma$ restricted to $E_k\times \Sigma$. Let us denote
it by $\cE_k=\E_k\stackrel{\uPhi_k}{\rightarrow}\E_k\times
\K$. Consider the kernel of $\uPhi_k$. Because $E_k$ parameterizes
nilpotent stable Higgs bundles with non-zero Higgs field
$\ker(\uPhi_k)$ is a line bundle over $E_k\times \Sigma$. Recall from 
Proposition 7.1 of \cite{hitchin1} that for $\compE\in N_k\subset E_k$ 
we have $\deg(\ker(\Phi))=1-k$. Since $E_k$ is smooth we have that 
$\deg(\ker(\Phi))=1-k$ for every $\compE\in E_k$.

The result follows.    
\end{proof}

\begin{remark} Clearly a completely analogous result holds for
$\tilde{N}$ with $\tN$, $\tilde{E}_k$ and $\tilde{N}_k$ instead of
$\N$, $E_k$ and $N_k$.
\end{remark}

\begin{notation} If $X$ is an irreducible locally closed
subvariety of a smooth algebraic variety $Y$ of
codimension d, 
then $\eta^Y_X\in H^{2d}(Y)$ denotes the cohomology class of $\overline{X}$ in
$Y$.  

If $X$ is an irreducible locally closed and relatively 
complete subvariety of $Y$ then
$\overline{\eta}^Y_X\in H_{cpt}^{2d}(Y)$ denotes the compactly supported
cohomology class of $\overline{X}$ in $Y$. 
\end{notation}

\begin{theorem} Let $k=1,..,g-1$. 
The degeneracy locus $D_k$ has the following
decomposition: $$D_k=\tN_k\cup\bigcup_{i=1}^k \tilde{E}^k_i,$$
where $\tN_k=D_k\cap \tN$, and $\tilde{E}_i^k\subset \tilde{E}_i$ are
those nilpotent stable Higgs bundles whose unique line bundle $L_\cE$ of 
Proposition~\ref{cone} has the property that $H^0(\Sigma,L_\cE\otimes
L_p^{k-1})\neq 0$.

Furthermore $\tilde{E}^k_k:=\{\cE\in\tilde{E}_k: L_\cE=L^{1-k}_p\}$
and hence 
\begin{eqnarray}
\eta^\tM_{\tilde{E}^k_k}\smallsetminus [\J]= \eta^\M_{E_k}
\in H^{6g-6}(\M) 
\label{splitE}
\end{eqnarray}
where $\eta^\tM_{\widetilde{E}^k_k}\smallsetminus [\J]$ means the
coefficient of $\eta^\J_{pt}$
in the decomposition of (\ref{splitM}).
\label{degeneracy}
\end{theorem} 

\begin{proof} Let $\cE=\compE$ be a stable Higgs bundle with 
$\Phi\neq 0$ and 
$\Hy^0(\Sigma,\cE\otimes L_p^{k-1})\neq 0$. It is easy to see that
this hypercohomology is the vector space of all morphisms from 
$\cE_0\otimes
L_p^{1-k}=L_p^{1-k}\stackrel{0}{\rightarrow}L_p^{1-k}\otimes K$ to $\cE$.
Consider a nonzero such morphism $f$. Consider $L$ the line subbundle
of $E$ generated by the image of $f$ of
Lemma~\ref{narasimhan-seshadri}. 
Clearly $L$ is killed by the Higgs field $\Phi$. This shows that
$\cE\in \tilde{N}$ and $L=L_\cE$. We also see that 
$\Hy^0(\Sigma,\cE\otimes L_p^{k-1})\cong H^0(L_\cE\otimes L_p^{k-1}$). The
first part of the statement follows. 

By the above argument we see that 
$\tilde{E}^k_k=\{\cE\in \tilde{E}_k: 
H^0(\Sigma, L_\cE \otimes L_p^{k-1})\neq 0\}$, however $L_\cE$ is of
degree $1-k$, thus  
$\tilde{E}^k_k=\{\cE\in\tilde{E}_k: L_\cE=L^{1-k}_p\}$, as claimed. 
This means that for every $\cE\in E_k$ there is a unique line bundle 
$L=L^{1-k}_p\otimes L^*_\cE$
such that $\cE\otimes L\in \tilde{E}^k_k$. This shows (\ref{splitE}).
\end{proof}

\begin{remark} By definition 
$\tilde{\N}_k=W^0_{2,2k-1}$ are non-Abelian Brill-Noether
loci as defined in  \cite{sundaram} (cf. \cite{teixidor}).
\end{remark}

\section{Proof of Theorem~\ref{main}}
\label{proof}

In this final section we prove Theorem~\ref{main}. 

\paragraph{\it Proof of Theorem~\ref{main}.} The proof proceeds by
showing that $\ch_0(\D_k)=4g-4$ then $c_{4g-3}(\D_k)=0$ and we
finish by using Porteous's theorem for $\D_k$. 

First we make some calculations.

\begin{lemma} The formal difference of coherent sheaves $\D_k$ has rank
$4g-4$, i.e. $$\ch_0(\D_k)=4g-4.$$ Moreover
\begin{eqnarray}
c(\D_k)=\left(1+\frac{\alpha_\M}{2}+\frac{\alpha^2_\M-\beta_\M}{4}\right)^{2g-2}
\label{chern}
\end{eqnarray} 
in the decomposition
(\ref{splitM}). 
\label{calc}
\end{lemma}

\begin{proof} It follows from the long exact sequence (\ref{long}) 
that 
$$\D_k=-{\pi_\tM}_!(\cE^k_\tM)={\pi_\tM}_!(\E^k_\tM\otimes \K)-
{\pi_\tM}_!(\E^k_\tM).$$ We can calculate the Chern character of the
right hand side by the \newline Grothendieck-Riemann-Roch theorem. This gives
$$\ch(\D_k)={\pi_\tM}_*\left(\ch(\E^k_\tM)(\ch
(\K)-1)\td(\Sigma)\right).$$
Now $\td(\Sigma)=1-(g-1)\sigma$ and $\ch(\K)=1+(2g-2)\sigma$.
Moreover ${\pi_\tM}_*$ maps a cohomology class in $a\in H^*(\tM)\otimes
H^*(\Sigma)$ of the form $$a=a_0\otimes 1 + \sum_{i=1}^{2g} a_1^i
\otimes e_i +
a_2\otimes \sigma$$ to the class $a_2\in H^*(\tM)$. The class $a_2$ is
denoted by $a_2=a\smallsetminus \sigma$, while the class $a_0$ is
denoted by $a_0=a\smallsetminus 1$. From this it follows that 
$$\ch(\D_k)=\left(\ch(\E^k_\tM)((2g-2)\sigma)(1-(g-1)\sigma)\right)
\smallsetminus\sigma=(2g-2)(\ch(\E^k_\tM)\smallsetminus 1).$$  

Observe that $\ch(\E^k_\tM)\smallsetminus 1 =\ch((\E^k_\tM)_p)\in
H^*(\tM)$, where $(\E^k_\tM)_p=\E^k_\tM\mid_{\tM \times \{p\}}$.
It follows from (\ref{c1}) and (\ref{c2}) 
that $c_1((\E^k_\tM)_p)=\alpha_\M$ and 
$c_2((\E^k_\tM)_p)=(\alpha_\M^2-\beta_\M)/4$. Hence the formal Chern
roots of $(\E^k_\tM)_p$ are $(\alpha_\M+\sqrt{\beta_\M})/2$ and 
$(\alpha_\M-\sqrt{\beta_\M})/2$. Thus 
\begin{eqnarray}
\ch((\E^k_\tM)_p)&=&\exp\left(\frac{\alpha_\M+\sqrt{\beta_\M}}{2}\right)
+\exp\left(\frac{\alpha_\M-\sqrt{\beta_\M}}{2}\right)\nonumber \\&=&
2e^{\alpha_\M/2} \cosh\left(\sqrt{\beta_\M}/2\right),\nonumber \end{eqnarray}
and hence 
$$\ch(\D_k)=(4g-4)e^{\alpha_\M/2}\cosh\left(\sqrt{\beta_\M}/2\right).$$
This shows that $\rank(\D_k)=\ch_0(\D_k)=4g-4$ and formal calculation
gives (\ref{chern}).
\end{proof}

(\ref{chern}) has the following immediate corollary:

\begin{corollary} $c_{4g-3}(\D_k)=0$.$\square$
\label{vanishing} 
\end{corollary}
 
To prove Theorem~\ref{main} we exhibit $g$ linearly independent elements 
$$r_0,r_1,..,r_{g-1}\in H^{6g-6}_{cpt}(\M)$$ for which $j_\M(r_i)=0$. 

To construct $r_k$ for $0<k<g$ consider the Zariski open subvarieties
$$\tM_k=\tM\setminus(\tN\bigcup_{i=1}^{k-1}\tilde{E}_i)$$ and 
$$\M_k=\M\setminus(\N\bigcup_{i=1}^{k-1}E_i)$$ of $\tM$ and $\M$ respectively. 
Restricting 
the sequence of Theorem~\ref{deg} to $\tM_k$ yields:
\begin{eqnarray}
0 &\longrightarrow& \R^0{\pi_\M}_*(\cE_{\tM}^k)\mid_{\tM_k}
\longrightarrow V\mid_{\tM_k}\stackrel{f\mid_{\tM_k}}{\longrightarrow} 
W\mid_{\tM_k} \nonumber \\ &\longrightarrow&
\R^1{\pi_\M}_*(\cE_{\tM}^k)\mid_{\tM_k}\longrightarrow
0.\label{sequ}
\end{eqnarray}

The degeneracy locus of $f\mid_{\tM_k}$ (where $f\mid_{\tM_k}$ fails
to be an injection) is $D_k\cap \tM_k$ which is $\tilde{E}^k_k$ from
Theorem~\ref{degeneracy}. This has codimension $4g-3$. Furthermore 
\begin{eqnarray}\rank(W)-\rank(V)&=&\rank\left(\R^1{\pi_\M}_*(\cE_{\tM}^k)\right)-
\rank\left(\R^0{\pi_\M}_*(\cE_{\tM}^k)\right)\nonumber \\ &=&\rank(\D_k)=4g-4 \nonumber\end{eqnarray}
by Lemma~\ref{calc}. Thus the degeneracy locus has the expected
codimension hence we are in the situation of Porteous's
theorem (cf. \cite{arbarello-et-al}), which gives: 
$$\eta_{\tilde{E}^k_k}^{\tM_k}=c_{4g-3}(W\mid_{\tM_k}-V\mid_{\tM_k})
\in H^{8g-6}(\tM_k).$$
The right hand side equals $c_{4g-3}(\D_k\mid_{\tM_k})$ by
(\ref{sequ}), which vanishes by Corollary~\ref{vanishing}. Moreover 
(\ref{splitE}) yields
$$\eta^{\tilde{\M}_k}_{\tilde{E}^k_k}
\smallsetminus [\J]= \eta^{\M_k}_{E_k}.$$ It follows that 
\begin{eqnarray}\eta^{\M_k}_{E_k}=0\in H^{6g-6}(\M_k).\label{0}\end{eqnarray}

From now on we work over $\M$. We show by induction 
on $i$ that there is
a formal linear combination $$r^i_k=\sum_{j=k-i}^k 
\lambda_j\cdot\left[\eta^{\M_{k-i}}_{E_j}\right]$$ of cohomology classes in 
$H^{6g-6}(\M_{k-i})$, such that $\lambda_k=1$ and 
the corresponding cohomology class 
$\sum_{j=k-i}^k\lambda_i\cdot\eta^{\M_{k-i}}_{E_j}$ is $0$ in
$H^{6g-6}(\M_{k-i})$. 

For $i=0$ the statement is just (\ref{0}). Suppose that there is 
such formal combination $r^i_k$. Consider the following bit of the
long exact sequence of the pair $\M_{k-i}\subset \M_{k-i-1}$:
$$ H^{6g-6}(\M_{k-i},\M_{k-i-1})\longrightarrow H^{6g-6}(\M_{k-i-1})
\longrightarrow H^{6g-6}(\M_{k-i}).$$
Because $\M_{k-i-1}\setminus \M_{k-i}=E_{k-i-1}$ is of real codimension
$6g-6$, the Thom isomorphism
transforms this sequence to:
\begin{eqnarray}
H^{0}(E_{k-i-1})\stackrel{\tau}{\longrightarrow} H^{6g-6}(\M_{k-i-1})
\stackrel{\rho}{\longrightarrow} H^{6g-6}(\M_{k-i}),\label{thom}
\end{eqnarray}
where $\tau$ is the Thom map and $\rho$ is restriction. Clearly 
$\rho(\eta^{\M_{k-i-1}}_{E_j})=\eta^{\M_{k-i}}_{E_j}$.
Thus $$\rho(\sum_{j=k-i}^k\lambda_j\cdot\eta^{\M_{k-i-1}}_{E_j})
=\sum_{j=k-i}^k\lambda_j\cdot\eta^{\M_{k-i}}_{E_j}=0.$$  
The exactness of (\ref{thom}) yields that the cohomology class 
$\sum_{j=k-i}^k \lambda_j\cdot\eta^{\M_{k-i-1}}_{E_j}$ is in the
image of $\tau$. Because $H^0(E_k)\cong \R$ there is a real number 
$-\lambda_{k-i-1}\in \R$ such that 
\begin{eqnarray}\tau(-\lambda_{k-i-1})=
\sum_{j=k-i}^k \lambda_j\cdot\eta^{\M_{k-i-1}}_{E_j}
\in H^{6g-6}(\M_{k-i-1})\label{akar}.
\end{eqnarray} 
However a well known property of the Thom map gives
$\tau(1)=\eta^{\M_{k-i-1}}_{E_{k-i-1}}$, thus from (\ref{akar}) the 
formal linear combination 
$$r^{i+1}_k=\sum_{j=k-i-1}^k 
\lambda_j\cdot\left[\eta^{\M_{k-i-1}}_{E_j}\right]$$ is
$0$, when considered as a class in $H^{6g-6}(\M_{k-i-1})$. This proves the 
existence of formal linear combinations $r^i_k$ for all $0\leq i\leq k-1$.

Using $r^{k-1}_k$ an identical argument gives the formal linear
combination $$r^\prime_k=\lambda\cdot\left[\eta^\M_\N\right]+
\sum_{j=1}^k \lambda_j\cdot\left[\eta^{\M}_{E_j}\right]$$ with the
property that $\lambda_k=1$ and $r^\prime_k$ when considered as an element of
$H^{6g-6}(\M)$ is $0$. Now the compactly supported cohomology class 
$$r_k=\lambda \cdot\overline{\eta}^\M_\N+
\sum_{j=1}^k \lambda_j\cdot \overline{\eta}^{\M}_{E_j}\in H^{6g-6}_{cpt}(\M)$$
has the property that $j_\M(r_k)=r^\prime_k=0$, where by abuse of
notation
$r^\prime_k$ denotes the cohomology class in $H^{6g-6}(\M)$
corresponding to the formal linear combination $r^\prime_k$.

We have found $g-1$ linearly independent compactly supported 
cohomology classes $r_1,..,r_{g-1}\in H^{6g-6}_{cpt}(\M)$. Clearly 
$\overline{\eta}^{\M}_\N$ is not in the span of $r_1,..,r_{g-1}$. Moreover
for each $0<i<g$ we have $$\int_\M \overline{\eta}^{\M}_\N \wedge
r_i=0$$ 
since
$j_\M(r_i)=0$. Furthermore   
$$\int_\M \overline{\eta}^{\M}_\N\wedge \overline{\eta}^{\M}_\N=\int_\N
c_{3g-3}(T^*_\N)= 0.$$ Thus $\overline{\eta}_\N^\M$ 
is perpendicular to  $r_1,..,r_{g-1}$
and $\overline{\eta}_\N^\M$, which constitutes a basis for $H^{6g-6}_{cpt}(\M)$, and
so  $j_\M(\overline{\eta}^{\M}_\N)=0$.

Putting our findings together: we have $g$ linearly independent
middle dimensional compactly supported classes $r_0=\overline{\eta}^{\M}_\N$
and $r_1,..,r_{g-1}$ 
in the kernel of the forgetful map 
$j_\M: H^{6g-6}_{cpt}(\M)\rightarrow H^{6g-6}(\M)$. 

Theorem~\ref{main} is finally proved.  {$\square$

\end{document}